\documentclass[11pt]{article}

\newcommand{\documentdate}{1 October 2015} 

\usepackage{a4wide,latexsym,mystyle,graphicx,varioref,color,amsmath,cite}
\DeclareGraphicsExtensions{.eps,.jpg}
\usepackage{epstopdf}

\title{Evaluation complexity bounds \\ for smooth constrained nonlinear optimization\\
  using scaled KKT conditions, high-order models and the criticality measure $\chi$}
\author{C. Cartis\thanks{
   Mathematical Institute, 
   Oxford University,
   Oxford OX2 6GG, Great Britain.
   Email: coralia.cartis@maths.ox.ac.uk.}
~ N. I. M. Gould\thanks{
   Numerical Analysis Group,
   Rutherford Appleton Laboratory,
   Chilton OX11 0QX, Great Britain.
   Email:  nick.gould@stfc.ac.uk.}
~and Ph. L. Toint\thanks{
   Namur Center for Complex Systems (naXys) and Department of Mathematics,
   University of Namur,                
   61, rue de Bruxelles, B-5000 Namur, Belgium.
   Email: philippe.toint@unamur.be.
   Ph. Toint would like to dedicace this paper to Tha\"{i}s, his first grand-daughter.}
}
\date{\documentdate}

\newcommand{\ass}[2]{\label{ass-#1}
                     \begin{list}{}{\setlength{\leftmargin}{2.5cm}}
                     \item \hspace{-2.5cm} \framebox[2.0cm]{\bf #1} \,\,#2
                     \end{list} }
\newcommand{\inneralg}{AR$p$CC}
\newcommand{\outeralg}{AR$p$GC}
\newcommand{\mycitetwo}[4]{[\citenum{#1},\,#2,\,\citenum{#3},\,#4]}
\newcommand{\ep}{\epsilon_{\mbox{\tiny P}}}
\newcommand{\ed}{\epsilon_{\mbox{\tiny D}}}
   
\begin{document}
\comment{

\begin{titlepage}

\includegraphics[height=4cm]{/Users/pht/Pictures/logos/UNamur.jpg}
\hspace*{4cm}
\includegraphics[height=3cm]{/Users/pht/Pictures/logos/naxys.jpg}
\vspace*{3cm}
\begin{center}
\begin{minipage}[c]{12cm}
\vfill
\begin{center}
{\sc  
Evaluation complexity bounds for \\smooth constrained nonlinear
optimization \\ using scaled KKT conditions, high-order models\\
and the $\chi$ criticality measure
}
\end{center}
\vfill
\centering{by C. Cartis, N. I. M. Gould and Ph. L. Toint}\\ 
\mbox{}
\vfill
\centering{\documentdate}\\
\vfill
\centering{\includegraphics[height=6cm]{smalleps.eps}}
\end{minipage}
\end{center}
\begin{center}
{\large
University of Namur, 61, rue de Bruxelles, B5000 Namur (Belgium)\\*[2ex]
{\tt http://www.unamur.be/sciences/naxys}}
\end{center}

\end{titlepage}
}

\maketitle

\begin{abstract}
\noindent
Evaluation complexity for convexly constrained optimization is considered and
it is shown first that the complexity bound of $O(\epsilon^{-3/2})$ proved by 
Cartis, Gould and Toint (IMAJNA 32(4) 2012, pp.1662-1695)
for computing an $\epsilon$-approximate first-order critical point can be
obtained under significantly weaker assumptions.  Moreover, the result is
generalized to the case where high-order derivatives are used, resulting in a
bound of $O(\epsilon^{-(p+1)/p})$ evaluations  whenever derivatives of order
$p$ are available.  It is also shown that the bound of
$O(\ep^{-1/2}\ed^{-3/2})$ evaluations ($\ep$ and $\ed$ being primal and dual
accuracy thresholds) suggested by Cartis, Gould and Toint (SINUM, 2015, to
appear) for the general nonconvex case involving both equality and inequality
constraints can be generalized to a bound of  $O(\ep^{-1/p}\ed^{-(p+1)/p})$
evaluations under similarly weakened assumptions.   This paper is variant of a
companion report (NTR-11-2015, University of Namur, Belgium) which uses a
different first-order criticality measure to obtain the same complexity bounds.
\end{abstract}

\vspace*{1cm}

\numsection{Introduction}

In \cite{CartGoulToin12b} and \cite{CartGoulToin13a}, we examined the worst-case
evaluation complexity of finding an $\epsilon$-approximate first-order
critical point for smooth nonlinear (possibly nonconvex) optimization
problems for a methods using both first and second derivatives of the
objective function.  The case where constraints are defined by a convex set
was considered in the first of these references while the general case (with
equality and inequality constraints) was discussed in the second. 

It was shown in \cite{CartGoulToin12b} that at most $O(\epsilon^{-3/2})$
evaluations of the objective function and its derivatives are needed to
compute such an approximate critical point.  This result, which is identical
in order to the best known result for the unconstrained case, comes at the
price of potentially restrictive technical assumptions: it was assumed that an
approximate first-order critical point of a cubic model subject to the
problem's constraints can be obtained for the subproblem solution in a
uniformly bounded number of descent steps that is independent of $\epsilon$,
that all iterates remains in a bounded set and that the gradient of the
objective function is also Lipschitz continuous (see \cite{CartGoulToin12b}
for details).  The analysis of \cite{CartGoulToin13a} then built on the result
of the convex case by first specializing it to convexly constrained nonlinear
least-squares and then using the resulting complexity bound in the context of
a two-phase algorithm for the problem involving general constraints.  If $\ep$
and $\ed$ are the primal and the dual criticality thresholds, respectively, it
was suggested that at most $O(\ep^{-1/2}\ed^{-3/2})$ evaluations of the
objective function and its derivatives are needed to compute an approximate
critical point in that case, where the Karush-Kuhn-Tucker (KKT) conditions are
scaled to take the size of the Lagrange multipliers into account. Because this
bound is based on that obtained for the convex case, it suffers from the same
limitations (not to mention an additional constraint on the relative sizes of
$\ep$ and $\ed$, see \cite{CartGoulToin13a}).

More recently, Birgin, Gardenghi, Mart\'{i}nez, Santos and Toint
\cite{BirgGardMartSantToin17}  provided a new regularization
algorithm for the unconstrained problem with two interesting features.
The first is that the model decrease condition used for the subproblem
solution is weaker than that used previously, and the second is that the
use of problem derivatives of order higher than two is allowed,
resulting in corresponding reductions in worst-case complexity.
In addition, the same authors also analyzed the worst-case evaluation
complexity of the general constrained optimization problem in
\cite{BirgGardMartSantToin16} also allowing for high-order derivatives and
models in a framework inspired by that of
\cite{CartGoulToin12e,CartGoulToin13a}.  At variance with the analysis of
these latter references, their analysis considers unscaled approximate
first-order critical points in the sense that such points satisfy the standard
unscaled KKT  conditions with accuracy $\ep$ and $\ed$. 

This paper is variant of a companion report \cite{CartGoulToin15a} which
uses a different first-order criticality measure to obtain the same complexity bounds. 
The first purpose of both papers is to explore the potential of the proposals
made in \cite{BirgGardMartSantToin17} to overcome the limitations of
\cite{CartGoulToin12b}  and to extend its scope by considering the use of
high-order derivatives and models. A second objective is to 
use the resulting worst-case bounds to establish strengthened evaluation
complexity bounds for the general nonlinearly constrained optimization problem
in the framework of scaled KKT conditions, thereby improving
\cite{CartGoulToin13a}.  This paper, like it companion, is thus organized in
two main sections, Section~\ref{convex-s} covering the convexly constrained
case and Section~\ref{general-s} that allowing general nonlinear constraints.
The results obtained are finally discussed in Section~\ref{concl-s}.
 
\numsection{Convex constraints}\label{convex-s}

The first problem we wish to solve is formally described as 
\beqn{prob}
\min_{x \in \calF} f(x)
\eeqn
where we assume that $f:\Re^n\longrightarrow\Re$ is $p$-times continuously 
differentiable, bounded from below, and has Lipschitz continuous 
$p$-th derivatives. For the $q$-the derivative of a function $h:\Re^n \to \Re$
to be Lipschitz continuous on the set $\calS \subseteq \Re^n$, we require that
i.e.\ there exists a constant $L_{h,q}\geq 0$ such that, for all $x,y \in
\calS$, 
\[
\|  \nabla_x^q h(x)  - \nabla_x^q h(y) \|_T \leq (p-1)! \, L_{h,p} \| x - y \|
\]
where $\|\cdot\|$ is the standard Euclidean norm on $\Re^n$ and $\|\cdot\|_T$
is recursively induced by this  norm on the space of $q$-th order tensors. We
also assume that the feasible set $\calF$ is closed, convex and
non-empty. Note that this formulation covers  unconstrained optimization
($\calF = \Re^n$), as well as standard inequality  
(and linear equality) constrained optimization in its different forms: 
the  set $\calF$ may be defined by simple bounds, and/or by polyhedral
or more general convex constraints. We are tacitly assuming here that the cost
of evaluating values and derivatives of the constraint functions possibly
involved in the definition of $\calF$ is negligible. 

The algorithm considered in this paper is iterative. Let $T_p(x_k,s)$ be the
$p$-th order Taylor-series approximation to $f(x_k+s)$ at some iterate $x_k
\in \Re^n$, and define the local regularized model at $x_k$ by 
\beqn{model}
m_k(x_k+s) \eqdef T_p(x_k,s)+\frac{\sigma_k}{p+1}\|s\|^{p+1},
\eeqn
where $\sigma_k> 0$ is the regularization parameter. Note that 
$m_k(x_k) = T_p(x_k,0) = f(x_k)$. The approach used in \cite{CartGoulToin12b}
(when $p=2$) seeks to define a new iterate $x_{k+1}$ from the preceding one by
computing an approximate solution of the subproblem 
\beqn{subprob}
\min_{x \in \calF} m_k(x_k+s)
\eeqn
using a modified version of the  Adaptive Regularization with Cubics (ARC)
method for unconstrained minimization. By contrast, we now examine the
possibility of modifying the AR$p$ algorithm of  \cite{BirgGardMartSantToin17}
with the aim of inheriting its interesting features. As in
\cite{CartGoulToin12b}, the modification involves a suitable continuous
first-order criticality measure for the constrained problem of minimizing a
given function $h: \Re^n \to \Re$ on $\calF$.  For an arbitrary $x \in \calF$,
this criticality measure is given by 
\beqn{chidef}
\chi_h(x) \eqdef  \left| \min_{x+d \in \calF, \| d \|_\chi \leq 1} \ip{\nabla_x h(x)}{d} \right|,
\eeqn
where $\ip{\cdot}{\cdot}$ denotes the Euclidean inner product and $\|\cdot\|_\chi$ 
is any fixed norm, possibly chosen to make the computation of $\chi_h(x)$ easier. Let $\kappa_n >0$ be the norm equivalence constant such that
\beqn{kndef}
\| v \| \leq \kappa_n \|v\|_\chi \tim{ for all} v \in \Re^n.
\eeqn
Observe that $\chi_h(x)$ depends on the geometry of $\calF$ only (and not on its possible parametrization using constraint functions) and that $x$ is a first-order critical point of problem \req{prob} if and only if $\chi_f(x) = 0$. Also note that $\chi_h(x) = \|\nabla_x h(x)\|$ whenever $\calF = \Re^n$ and $\|\cdot\|_\chi = \|\cdot\|$.

We now describe our algorithm as the \inneralg\ algorithm
(AR$p$ for Convex Constraints) \vpageref{\inneralg}.

\vspace*{-0.2cm}
\algo{\inneralg}{Adaptive Regularization using 
$p$-th order models  for convex constraints (\inneralg)}{
A starting point $x_{-1}$, an initial and a minimal regularization parameter 
$\sigma_0\geq \sigma_{\min}>0$,  algorithmic parameters $\theta > 0$,
$\gamma_3\geq \gamma_2 > 1 > \gamma_1 >0$ and 
$1>\eta_2\geq \eta_1>0$, are given, as well as an accuracy threshold $\epsilon\in (0,1]$.
Compute $x_0 = P_\calF[x_{-1}]$, the projection of $x_{-1}$ onto $\calF$,
and evaluate $f(x_0)$ and $\nabla_x f(x_0)$.

For $k=0, 1, \ldots$, do:
\begin{enumerate}
\item Evaluate $\nabla_x f(x_k)$.
If 
\beqn{terminateg}
\chi_f(x_k)\leq \epsilon,
\eeqn
terminate with $x_\epsilon = x_k$. Otherwise compute derivatives of $f$ of order 2 to $p$ at $x_k$.
\item Compute a step $s_k$ by approximately minimizing $m_k(x_k+s)$
over $s \in \calF$ so that
\beqn{sub-term-1}
x_k+s_k \in \calF,
\eeqn
\beqn{sub-term-2}
m_k(x_k+s_k)<m_k(x_k)
\eeqn
and
\beqn{sub-term-3}
\chi_{m_k}(x_k+s_k)  \leq \theta \, \|s_k\|^p.
\eeqn
\item Compute $f(x_k+s_k)$ and 
\beqn{rhodef}
\rho_k = \frac{f(x_k)-f(x_k+s_k)}{T_p(x_k,0)-T_p(x_k,s_k)}.
\eeqn
If $\rho_k \geq \eta_1$, set $x_{k+1} = x_k+s_k$. Otherwise set $x_{k+1}=x_k$.
\item Set
\beqn{vsu}
\sigma_{k+1} \in \left\{ \arr{cll}{
[\max(\sigma_{\min}, \gamma_1\sigma_k)\sigma_k ] & \tim{if} \rho_k > \eta_2 &\hspace*{-0.3cm}\tim{[very successful iteration]} \\
\;[ \sigma_k, \gamma_2 \sigma_k ] & \tim{if} \eta_1 \leq \rho_k \leq \eta_2& 
\hspace*{-0.3cm}\tim{[successful iteration]}\\
\; [ \gamma_2 \sigma_k, \gamma_3 \sigma_k ] &  \tim{otherwise.} & \hspace*{-0.3
cm}\tim{[unsuccessful iteration],}
}\right.
\eeqn
and go to step 2 if $\rho_k < \eta_1$.
\end{enumerate}
}

  We first state a useful property of the \inneralg\ algorithm, which ensures
  that a fixed fraction of the iterations $1, 2, \ldots, k$ must be either
  successful or very successful. 

\llem{unsuccbounbd}{ 
\mycitetwo{BirgGardMartSantToin17}{Lem.2.4}{CartGoulToin12e}{Thm.2.2}.
Assume that, for some $\sigma_{\max} > 0$,  $\sigma_j \leq \sigma_{\max}$
for all $0\leq j \leq k$. Then the \inneralg\ algorithm ensures that
\beqn{kbound}
k \leq \kappa_u |\calS_k|, \tim{where} \kappa_u \eqdef
\left( 1 + \frac{|\log \gamma_1|}{\log \gamma_2}
 \right) + \frac{1}{\log \gamma_2} 
 \log\left(\frac{\sigma_{\max}}{\sigma_0}\right),
 \eeqn
where $\calS_k$ is the number of successful  and very successful iterations,
in the sense of \req{vsu},  up to iteration $k$.}

We start our worst-case analysis by formalizing our assumptions, using
$\calL(x_0) = \{ x \in \calF \mid f(x) \leq f(x_0) \}$

\ass{AS.1}{
The objective function $f$ is $p$ times continuously differentiable on an open
set containing $\calL(x_0)$. 
}

\ass{AS.2}{ The $p$-th derivative of $f$ is Lipschitz continuous on $\calL(x_0)$.
}

\ass{AS.3}{ The feasible set $\calF$ is closed, convex and non-empty.}

The \inneralg\ algorithm is required to start from a feasible $x_0 \in \calF$,
which, together with the fact that the subproblem solution in Step~2 involves
minimization over $\calF$, leads to AS.3.

We now recall some simple results whose proof can be found in
\cite{BirgGardMartSantToin17} in the context of the original AR$p$ algorithm. 

\llem{ARplemma}{
Suppose that AS.1 and AS.2 hold.  Then, for each $k \geq 0$,


\begin{enumerate}
\item[(i)] 
\beqn{errf}
f(x_k+s_k) \leq T_p(x_k,s_k) + \frac{L_{f,p}}{p} \|s_k\|^{p+1}
\eeqn
and 
\beqn{errg}
\|\nabla_x f(x_k+s_k) - \nabla_s T(x_k,s_k) \| \leq L _{f,p} \|s_k\|^p;
\eeqn
\item[(ii)] 
\beqn{Dphi}
T_p(x_k,0) - T_p(x_k,s_k) \geq \frac{\sigma_k}{p+1} \|s_k\|^{p+1};
\eeqn
\item[(iii)]
\beqn{sigmaupper}
\sigma_k \leq \sigma_{\max}
\eqdef \max\left[ \sigma_0, \frac{\gamma_3 L_{f,p}(p+1)}{p \,(1-\eta_2)} \right].
\eeqn
\end{enumerate}
}

\proof{
See \cite{BirgGardMartSantToin17} for the proofs  of \req{errf} and
\req{errg}, which crucially depend on AS.1 and AS.2 being valid on the segment
$[x_k, x_k+s_k]$. Observe also that \req{model} and \req{sub-term-2} ensure
\req{Dphi}. 
Assume now that
\beqn{siglarge}
\sigma_k \geq \frac{L_{f,p}(p+1)}{p \,(1-\eta_2)}.
\eeqn
Using \req{errf} and \req{Dphi}, we may then deduce that
\[
|\rho_k - 1|
\leq \frac{|f(x_k+s_k) - T_p(x_k,s_k)|}{|T_p(x_k,0)-T_p(x_k,s_k)|}
\leq \frac{L_{f,p}(p+1)}{p \,\sigma_k} 
\leq 1-\eta_2
\]
and thus that $\rho_k \geq \eta_2$. Then iteration $k$ is very successful in
that $\rho_k \geq \eta_2$ and $\sigma_{k+1}\leq \sigma_k$.  As a consequence,
the mechanism of the algorithm ensures that \req{sigmaupper} holds. 
}

We now prove that, at successful iterations, the step at iteration $k$ must be
bounded below by a multiple of the $p$-th root of the criticality measure at
iteration $k+1$. 

\llem{schi+-lemma}{
Suppose that AS.1--AS.3  hold. Then 
\beqn{schi+}
\|s_k\| \geq  \left[ \frac{\chi_f(x_{k+1})}{2\kappa_n(L_{f,p}+\theta+\sigma_{\max})}\right]^{\frac{1}{p}}
\tim{ for all } k \in \calS.
\eeqn
}

\proof{
Since $k\in \calS$ and by definition of the trial point, we have that $x_{k+1}= x_k+s_k$.
Observe now that \req{errg} and  \req{sigmaupper} imply that
\beqn{errgradm}
\|\nabla f(x_{k+1}) - \nabla_x m_k(x_{k+1})\|
 \leq L_{f,p} \|s_k\|^p + \sigma_k \|s_k\|^p
 \leq (L_{f,p} + \sigma_{\max}) \|s_k\|^p,
\eeqn
and also that
\beqn{con2-ARC2CC-sl-1}
\begin{array}{lcl}
\chi_f(x_{k+1}) 
& \eqdef & |\ip{\nabla_x f(x_{k+1})}{d_{k+1}}| \\*[1ex]
&  \leq  & |\ip{\nabla_x f(x_{k+1}) - \nabla_s m_k(x_{k+1})}{d_{k+1}}| 
           + |\ip{\nabla_s m_k(x_{k+1})}{d_{k+1}}|,
\end{array}
\eeqn
where the first equality defines the vector $d_{k+1}$ with
\beqn{dplusbound}
\|d_{k+1}\|_\chi \leq 1.
\eeqn

Assume now, for the purpose of deriving a contradiction, that \req{schi+}
fails at iteration $k \in \calS$. 
Using the Cauchy-Schwarz inequality, \req{kndef}, \req{dplusbound},
\req{errgradm}, the failure of \req{schi+}  and the first
part of \req{con2-ARC2CC-sl-1} successively, we then obtain that
\[
\begin{array}{lcl}
\lefteqn{\ip{\nabla_s m_k(x_{k+1})}{d_{k+1}} - \ip{\nabla_x f(x_{k+1})}{d_{k+1}}}\\*[1ex]
\hspace*{4cm} &\leq &|\ip{\nabla_x f(x_{k+1})}{d_{k+1}} - \ip{\nabla_s m_k(x_{k+1})}{d_{k+1}}| \\*[1ex]
&\leq &\|\nabla_x f(x_{k+1}) - \nabla_s m_k(x_k+s_l)\| \, \|d_{k+1}\| \\*[1ex]
&\leq &\kappa_n ( L_{p,f}+ \sigma_{\max}) \|s_k\|^p \\*[1ex]
&\leq &\kappa_n ( L_{p,f}+ \theta+\sigma_{\max}) \|s_k\|^p \\*[1ex]
&\leq &\half \chi_f(x_{k+1}) \\*[1ex]
&  =  & - \half \ip{\nabla_x f(x_{k+1})}{d_{k+1}},
\end{array}
\]
which in turn ensures that
\[
\ip{\nabla_s m_k(x_{k+1})}{d_{k+1}} \leq \half \ip{\nabla_x f(x_{k+1})}{d_{k+1}} < 0.
\]
Moreover, $x_{k+1}+d_{k+1} \in \calF$ by definition of $\chi_f(x_{k+1})$,
and hence, using \req{dplusbound},
\beqn{nmdleqchi}
|\ip{\nabla_s m_k(x_{k+1})}{d_{k+1}}| \leq \chi_{m_k}(x_{k+1}).
\eeqn
We may then substitute this inequality in \req{con2-ARC2CC-sl-1} and use the
Cauchy-Schwarz inequality, \req{kndef} and \req{dplusbound} again to deduce
that 
\beqn{con2-ARC2CC-sl-1b}
\chi_f(x_{k+1}) 
\leq  \|\nabla_x f(x_{k+1}) - \nabla_s m_k(x_{k+1})\| + \chi_{m_k}(x_{k+1})
\leq  \kappa_n ( L_p + \alpha + \sigma_{\max}) \|s_k\|^p
\eeqn
where the last inequality results from \req{errgradm}, the identity $x_{k+1} =
x_k+s_k$ and \req{sub-term-3}.  But this contradicts our assumption that
\req{schi+} fails. Hence \req{schi+} must hold. 
}

We now consolidate the previous results by deriving a lower bound on the
objective function decrease at successful iterations. 
 
\llem{decrsucc}{
Suppose that AS.1--AS.3 hold.  Then, if iteration $k$ is successful, 
\[
f(x_k)-f(x_{k+1})  \geq \frac{1}{\kappa_s^f} \,\chi_f(x_{k+1})^{\frac{p+1}{p}}
\]
where
\beqn{kappasdef}
\kappa_s^f 
\eqdef \frac{p+1}{\eta_1 \sigma_{\min}} \Big[2\kappa_n(L_{f,p}+\theta+\sigma_{\max})\Big]^{\frac{p+1}{p}}.
\eeqn
}

\proof{
  If iteration $k$ is successful, we have, using \req{rhodef},
  \req{Dphi}, \req{vsu}, \req{schi+} and \req{sigmaupper} successively, that
\[
\begin{array}{lcl}
f(x_k)-f(x_{k+1}) 
& \geq & \eta_1 [\, T_p(x_k,0)-T_p(x_k,s_k)\, ] \\*[3ex]
&\geq &\bigfrac{\eta_1 \sigma_{\min}}{p+1} \;\|s_k\|^{p+1} \\*[3ex]
&\geq &\bigfrac{\eta_1 \sigma_{\min}}{(p+1)[2\kappa_n(L_{f,p}+\theta + \sigma_{\max})]^{\frac{p+1}{p}} }\;
             \chi_f(x_{k+1})^{\frac{p+1}{p}}.
\end{array}
\]
} 

\noindent
It is important to note that the validity of this lemma does not depend on the
history of the algorithm, but is only conditional to the smoothness assumption
on the objective function holding along the step from $x_k$ to $x_{k+1}$.  We
will make use of that observation in Section~\ref{general-s}. 

Our worst-case evaluation complexity results can now be proved by combining
this last result with the fact that $\chi_f(x_k)$ cannot be smaller than
$\epsilon$ before termination. 

\lthm{final_theorem}{Suppose that AS.1--AS.2 hold and let $f_{\rm low}$ be a
  lower bound on $f$ on $\calF$.  Then, given $\epsilon> 0$, the
  \inneralg\ algorithm applied on problem \req{prob} needs at most 
\[
\left \lfloor \kappa_s^f
                 \;\frac{f(x_0)- f_{\rm low}}{\epsilon^{\frac{p+1}{p}}} \right \rfloor
\]
successful iterations (each involving one evaluation of $f$ and its $p$ first derivatives) 
and at most
\[
\kappa_u \left \lfloor \kappa_s^f
                 \;\frac{f(x_0)- f_{\rm low}}{\epsilon^{\frac{p+1}{p}}} \right \rfloor
\]
iterations in total to produce an iterate $x_\epsilon$ such that
$\chi_f(x_\epsilon) \leq \epsilon$, where $\kappa_u$ is given by \req{kbound}
with $\sigma_{\max}$ defined by \req{sigmaupper}. 
}

\proof{
At each successful iteration, we have, using Lemma~\ref{decrsucc}, that
\[
f(x_k)-f(x_{k+1}) 
\geq  (\kappa_s^f)^{-1}\chi_f(x_{k+1})^{\frac{p+1}{p}}
\geq  (\kappa_s^f)^{-1}\epsilon^{\frac{p+1}{p}},
\]
where we used the fact that $\chi_f(x_{k+1}) \geq \epsilon$ before termination
to deduce the last inequality. Thus we deduce that, as long as termination
does not occur, 
\[
f(x_0) - f(x_{k+1})
= \sum_{j\in \calS_k} [f(x_j)-f(x_j+s_j) ]
\geq \bigfrac{|\calS_k|  }{\kappa_s^f}\;\epsilon^{\frac{p+1}{p}},
\]
from which the desired bound on the number of successful iterations
follows. Lemma~\ref{unsuccbounbd} is then invoked to compute the upper bound
on the total number of iterations. 
}

\numsection{The general constrained case} \label{general-s}

We now consider the general smooth constrained problem in the form
\beqn{genprob}
\min_{x \in \calF} f(x) 
\tim{ subject to }
c(x) = 0
\eeqn
where $c:\Re^n \to \Re^m$ is sufficiently smooth and $f$ and $\calF$ are as above.
Note that this formulation covers the general problem involving both equality
and inequality constraints, the latter being handled using slack variables and
the inclusion of the associated simple bounds in the definition of $\calF$.
In order to revise our smoothness assumptions, we first define, for some
parameter $\beta >0$,  the  neighbourhood of the feasible set given  by
\[
\calC_\beta= \{ x \in \calF \mid \|c(x)\| \leq \beta \},
\]
where our revised assumptions on the objective function have to hold (we
continue to assume AS.3). 

\ass{AS.4}{The objective function $f$ is $p$ times continuously differentiable
  on an open set containing $\calC_\beta$. 
}

\ass{AS.5}{All derivatives of $f$ of order 1 to $p$ are uniformly bounded and
  Lipschitz continuous in $\calC_\beta$. 
}

\ass{AS.6}{For each $i = 1, \ldots,m$, the constraint function $c_i$ is $p$
  times continuously differentiable on an open set containing $\calF$. 
}

\ass{AS.7}{All derivatives of order 1 to $p$ of each $c_i$ ($i=1,\ldots,m$)
  are uniformly bounded and Lipschitz continuous in $\calF$. 
}

\ass{AS.8}{There exists constants $f_{\rm low} \leq f_{\rm up}$ such that
$f(x) \in [f_{\rm low}, f_{\rm up}]$ for all $x \in \calC_\beta$.
}

\noindent
Note that AS.3, AS.5 and AS.7 allow us to apply the \inneralg\ algorithm to
the problem 
\beqn{feas-prob}
\min_{x\in \calF} \half \|c(x)\|^2.
\eeqn 
for any $\epsilon \leq \beta$.
If an approximately feasible point is found, then, because of AS.3--AS.7,  the
same \inneralg\ may then be applied to approximately solve the problem 
\beqn{phase2-problem}
\min_{x \in \calF} \mu(x,t_k) 
\eqdef \half \|r(x,t_k)\|^2  
\eqdef \half \left\|\left( \begin{array}{c}
 c(x) \\ f(x) - t_k \end{array} \right) \right\|^2
\eeqn
for some monotonically decreasing sequence of ``targets'' $t_k$ ($k=1,\ldots$).
This suggests that we might solve the problem \req{phase2-problem} using a
two-phase algorithm much in the spirit of that proposed by Cartis \textit{et
  al.} \cite{CartGoulToin11b,CartGoulToin12e,CartGoulToin13a} and Birgin
\textit{et al.} \cite{BirgGardMartSantToin17}. It is described
\vpageref{\outeralg}.  

\algo{\outeralg}{Adaptive Regularization using 
$p$-th order models  for general constraints (\outeralg)}{
  A constant $\beta$ defining $\calC_\beta$, a starting point $x_{-1}$,
  a minimum regularization parameter $\sigma_{\min} > 0$, an initial
  regularization parameter $\sigma_0\geq \sigma_{\min}$ are given,  as well
as a constant $\delta > 1$. The primal and dual tolerances   
\[
0 < \ep \leq \min\left[\beta, \left(\frac{\delta-1}{\delta}\right)^p, 1 \right]
\tim{ and }
0 < \ed < 1
\]
are also given.\\
\begin{description}
\item[Phase 1: ] \mbox{}\\
Starting from $x_0=P_\calF(x_{-1})$, apply the \inneralg\  algorithm to minimize
$\half \|c(x)\|^2$ subject to $x \in \calF$ until a point $x_1 \in \calF$ is found such that
\beqn{phase1-term}
\|c(x_1)\| \leq \ep -\ep^{\frac{p+1}{p}} \tim{ or } \chi_{\half\|c\|^2}(x_1) \leq \ed \|c(x_1)\|.
\eeqn
If $\|c(x_1)\| > \ep -\ep^{\frac{p+1}{p}} $, terminate with $x_\epsilon = x_1$.

\item [Phase 2: ] \mbox{}
\begin{enumerate}
\item Set $t_1 = f(x_1)-\sqrt{\ep^2-\|c(x_1)\|^2}$.

\item For $k=1, 2, \ldots$, do:
\begin{enumerate}
\item Starting from $x_k$, apply the \inneralg\ algorithm
   to minimize $\mu(x,t_k)$ as a function of $x \in \calF$ until 
   an iterate $x_{k+1} \in \calF$ is found such that 
   \beqn{ph2term}
   \|r(x_{k+1},t_k)\| \leq  \ep-\ep^{\frac{p+1}{p}}
   \tim{ or }
    f(x_{k+1}) < t_k
    \tim{ or }
   \chi_\mu(x_{k+1},t_k) \leq \ep\ed 
  \eeqn
\item 
\begin{enumerate}
\item If $\|r(x_{k+1},t_k)\| <  \ep-\ep^{\frac{p+1}{p}}$,
define $t_{k+1}$ according to 
\beqn{tk-update}
t_{k+1}=f(x_{k+1})-\sqrt{\ep^2-\|c(x_{k+1})\|^2}.
\eeqn
and terminate with $(x_\epsilon,t_\epsilon) = (x_{k+1},t_{k+1})$ if 
$\chi_\mu(x_{k+1},t_{k+1}) \leq \ep \ed$.\\
\item If  $\|r(x_{k+1},t_k)\| \geq  \ep-\ep^{\frac{p+1}{p}}$ and  $f(x_{k+1}) < t_k$, 
define $t_{k+1}$ according to
\beqn{tk-swap}
t_{k+1} = 2 f(x_{k+1})-t_k
\eeqn
and terminate with $(x_\epsilon,t_\epsilon) = (x_{k+1},t_{k+1})$ if 
$\chi_\mu(x_{k+1},t_{k+1}) \leq \ep \ed$.\\
\item If $\|r(x_{k+1},t_k)\| \geq  \ep-\ep^{\frac{p+1}{p}}$ and $f(x_{k+1}) \geq  t_k$,
terminate with $(x_\epsilon,t_\epsilon) = (x_{k+1},t_k)$
\end{enumerate}
\end{enumerate}
\end{enumerate}
\end{description}
}

Observe that the recomputations of $\chi_\mu(x_{k+1},t_{k+1})$ in Step~2.(b)
do not require re-evaluating $f(x_{k+1})$ or $c(x_{k+1})$ or any of their
derivatives. 

We now start our analysis by examining the complexity of Phase~1.

\llem{ph1-compl}{
Suppose that AS.3, AS.4 and AS.6 hold.  Then Phase 1 of the
\outeralg\ algorithm terminates after at most 
\[
\left\lfloor
\kap{CC}^c \|c(x_0)\| \, \ep^{-\frac{1}{p}} \, \ed^{-\frac{p+1}{p}}
\right\rfloor
\]
evaluations of $c$ and its derivatives, where $\kap{CC}^c
\eqdef \half \kappa_u \kappa_s^{\half\|c\|^2}\delta^{\frac{1}{p}}$. 
}

\proof{
Let us index the iteration of the \inneralg\ algorithm applied on problem \req{feas-prob} by $j$.
Assume that iteration $j$ is successful and that 
\beqn{biggish}
\|c(x_j)\| > \ep -\ep^{\frac{p+1}{p}} = \ep ( 1 - \ep^{1/p} )
\geq \delta^{-1} \ep,
\eeqn
\clearpage
where the last inequality follows from the bound on $\ep$ as a function of $\delta$.
Then, using the decreasing nature of the sequence $\{\|c(x_j)\|\}$,
Lemma~\ref{decrsucc} and the second part of \req{phase1-term}, we obtain that 
\[
 (\|c(x_j)\| - \|c(x_{j+1})\|)\,\|c(x_j)\|
\geq  \half \|c(x_j)\|^2 -  \half \|c(x_{j+1})\|^2
\geq  \left(\kappa_s^{\half\|c\|^2}\right)^{-1}  (\ed \|c(x_j)\|)^{\frac{p+1}{p}}
\]
and thus that
\[
\|c(x_j)\| - \|c(x_{j+1})\|
\geq  \left(\kappa_s^{\half\|c\|^2}\right)^{-1} \,\|c(x_j)\|^{\frac{1}{p}} \,\ed^{\frac{p+1}{p}}
\geq  \left(\kappa_s^{\half\|c\|^2}\right)^{-1} \delta^{-\frac{1}{p}} \,\ep^{\frac{1}{p}}\,\ed^{\frac{p+1}{p}}
\]
where we have used \req{biggish} to derive the last inequality. As in
Theorem~\ref{final_theorem}, we then deduce that the number of successful
iterations required for the \inneralg\ algorithm to produce a point $x_1$
satisfying \req{phase1-term} is bounded above by 
\[
\half \kappa_s^{\half\|c\|^2} \, \delta^{-\frac{1}{p}} \|c(x_0)\| \,
\ep^{-\frac{1}{p}}\,\ed^{-\frac{p+1}{p}}.
\]
The desired conclusion the follows by using Lemma~\ref{unsuccbounbd}.
} 


\noindent
We now partition the Phase~2 outer iterations into two subsets whose indexes
are given by 
\beqn{K-pos}
\calK_+ 
\eqdef \{ k \geq 0 \mid   \mbox{\req{tk-update}  holds} \}
\tim{ and }
\calK_- 
\eqdef  \{ k \geq 0 \mid  \mbox{\req{tk-swap} holds}  \}.
\eeqn
This partition allows us to prove the following technical results.

\llem{tech-l}{
The sequence  $\{t_k\}$ is monotonically decreasing.  Moreover, in every Phase
2 iteration of the \outeralg\ algorithm of index $k \geq 1$, we have that 
\beqn{tkltfk}
f(x_k) - t_k \geq 0,
\eeqn
\beqn{nreseps}
\|r(x_{k+1}, t_{k+1})\| = \ep  \tim{ for } k \in \calK_+,
\eeqn
\beqn{r-swap}
\| r(x_{k+1},t_{k+1})\| = \| r(x_{k+1},t_k)\| \leq \ep \tim{ for } k \in \calK_-,
\eeqn
\beqn{approxfeas}
\|c(x_k)\| \leq \ep \tim{and} f(x_k) - t_k \leq \ep,
\eeqn 
\beqn{tkdecr}
t_k - t_{k+1} \geq \ep^{\frac{p+1}{p}} \tim{ for } k \in \calK_+.
\eeqn
Moreover, if AS.8 holds, then, for $k \geq 1$,
\beqn{tk-bounded}
t_k \in [f_{\rm low}-1,f_{\rm up}]. 
\eeqn
Finally, at termination of the \outeralg\ algorithm,
\beqn{conds-at-termination}
\|r(x_\epsilon,t_\epsilon)\| \geq  \ep-\ep^{\frac{p+1}{p}}
\tim{ and }
f(x_\epsilon) \geq t_\epsilon
\tim{ and }
\chi_\mu(x_\epsilon,t_\epsilon) \leq \ep \ed.
\eeqn
}

\proof{
The inequality \req{tkltfk} follows from \req{tk-update} for $k-1 \in \calK_+$
and from \req{tk-swap} for $k-1 \in \calK_-$.  \req{nreseps} is also deduced
from \req{tk-update} while  \req{tk-swap} implies the equality in
\req{r-swap}, the inequality in that statement resulting from the decreasing
nature of $\|r(x,t_k)\|$ during inner iterations in Step~2.(a) of the
\outeralg\ algorithm. The inequalities \req{approxfeas} then follow from
\req{tkltfk},  \req{nreseps} and \req{r-swap}. 
We now prove \req{tkdecr}, which only occurs when
$\|r(x_{k+1},t_k)\| \leq \ep-\ep^{\frac{p+1}{p}}$, that is when
\beqn{res1}
(f(x_{k+1}) - t_k)^2 + \|c(x_{k+1})\|^2 \leq \left( \ep-\ep^{\frac{p+1}{p}}\right)^2.
\eeqn
From \req{tk-update}, we then have that
\beqn{dtk}
t_k-t_{k+1}= -(f(x_{k+1})-t_k)+\sqrt{\|r(x_k,t_k)\|^2-\|c(x_{k+1})\|^2}.
\eeqn
Now taking into account that the global minimum  of the problem 
\[
\min_{(f,c)\in\smallRe^2} \psi(f,c)\eqdef -f+\sqrt{\ep^2 - c^2} 
\tim{subject to} f^2+c^2 \leq \omega^2,
\]
for $\omega \in [0, \ep]$ is attained at $(f_*,c_*)=(\omega, 0)$ and it is given by
$\psi(f_*,c_*)=\ep-\omega$ (see \cite[Lemma 5.2]{CartGoulToin13a}), we obtain
from \req{res1} and \req{dtk} (setting $\omega = \ep-\ep^{\frac{p+1}{p}}$)
that 
\[
t_k-t_{k+1} \geq \ep - \omega = \ep^{\frac{p+1}{p}},
\]
for $k \in \calK_+$, which is \req{tkdecr}.
Note that, if $k \in \calK_-$, then we must have that $t_k > f(x_{k+1})$ and
thus \req{tk-swap} ensures that $t_{k+1} < t_k$. This observation and
\req{tkdecr} then allow us to conclude that the sequence $\{t_k\}$ is
monotonically decreasing. 

The inclusion \req{tk-bounded} is deduced from Step~1 of Phase~2 of the
\outeralg\ algorithm, the decreasing nature of the sequence $\{t_k\}$,
\req{tkltfk}, \req{approxfeas} and AS.8. 

In order to prove \req{conds-at-termination}, we need to consider, in turn,
each of the three possible cases where termination occurs in Step~2.(b).  In
the first case (i), $\|r(x_{k+1},t_k)\|$ is small (in the sense that the first
inequality of \req{ph2term} holds) and \req{tk-update} is then used, implying
that \req{nreseps} holds and that $f(x_{k+1}) > t_{k+1}$. If termination
occurs because $\chi(x_{k+1},t_{k+1}) \leq \ep \ed$, then
\req{conds-at-termination} clearly holds at $(x_{k+1},t_{k+1})$.  In the
second case (ii), $\|r(x_{k+1},t_k)\|$ is large (the first inequality of
\req{ph2term} fails), but $f(x_{k+1}) < t_k$, and $t_{k+1}$ is then defined by
\req{tk-swap}, ensuring that $f(x_{k+1}) > t_{k+1}$ and, because of
\req{r-swap}, that $\|r(x_{k+1},t_{k+1})\|$ is also large.  As before
\req{conds-at-termination}  holds at $(x_{k+1},t_{k+1})$ if termination occurs
because $\chi(x_{k+1},t_{k+1}) \leq \ep \ed$. The third case (iii) is when
$\|r(x_{k+1},t_k)\|$ is sufficiently large and $f(x_{k+1}) \geq t_k$.  But
\req{ph2term} then guarantees that $\chi(x_{k+1},t_k) \leq \ep \ed$, and the
inequalities \req{conds-at-termination} are again satisfied at
$(x_{k+1},t_k)$. 
} 

\noindent
Using the results of this lemma allows us to bound the number of outer iterations in $\calK_+$. 

\llem{number-outer}{
Suppose that AS.3, AS.4, AS.6  and AS.8 hold. Then 
\[
| \calK_+ | \leq [f_{\rm up}-f_{\rm low} + 1]\,\ep^{-\frac{p+1}{p}}.
\]
}

\proof{
We know from Lemma~\ref{tech-l} that $t_k$ decreases monotonically with, by \req{tkdecr}, a decrease of at least $\ep^{\frac{p+1}{p}}$ for $k \in \calK_+$.  Hence the desired conclusion follows from \req{tk-bounded}.
} 

\noindent
We now state a very useful consequence of \req{tk-bounded}, which is of interest for the analysis of the inner iterations.

\llem{uniform-Lsig}{
Suppose that AS.4-AS.8 hold. Then there exists a constant $L_{\mu,p} \geq 0$
such that the $p$-th derivative of $\mu(x,t_k)$ with respect to $x$ is
Lipschitz continuous with Lipschitz constant $L_{\mu,p}$ for all values of
$t_k$ computed by the \outeralg\ algorithm.  Furthermore, there exists a
constant $\sigma_{\mu,\max} > \sigma_{\min}$ such that all regularization
parameters arising in the \inneralg\ algorithm within Step~2.(a) of the
\outeralg\ algorithm are bounded above by $\sigma_{\mu,\max}$. 
}

\proof{
Because of AS.4--AS.7, we obtain that, for any 
$t \in [f_{\rm low}-1,f_{\rm up}]$,  
$\mu(x,t)$ is $p$ times continuously differentiable in $\calC_\beta$ and its
$p$-th derivative is Lipschitz continuous.  Moreover, since $\mu(x,t)$ and its
derivatives depend continuously on $t \in [f_{\rm low}-1,f_{\rm up}]$,  we may
deduce the existence of $L_{\mu,p}$, which is an upper bound on the $p$-th
derivative Lipschitz constant associated with each $t$. This proves the first
part of the proposition. The second is then derived by introducing $L_{\mu,p}$
in \req{sigmaupper} as specified by Lemma~\ref{ARplemma} to obtain
$\sigma_{\mu,\max}$. 
} 

\noindent
The main consequence of this result is that we may apply the
\inneralg\ algorithm to the minimization of $\mu(x,t_k)$ in Step~2.(a) of the
\outeralg\ algorithm and use all the properties of the former (as derived in
the previous section) using problem constants valid for every possible $t_k$,
because of \req{tk-bounded}. 

Consider now $x_k$ for $k \in \calK_+$ and  denote by $x_{k+\ell(k)}$ the next
iterate such that $k+\ell(k) \in \calK_+$. Two cases are then possible
(assuming termination does not occur at $x_{k+\ell(k)}$): either a single pass
in Step~2.(a) of the \outeralg\ algorithm is sufficient to obtain
$x_{k+\ell(k)}$ ($\ell(k) = 1$) or  two or more passes are necessary, with
iterations $k, \ldots, k+\ell(k)-1$ belonging to $\calK_-$.  Assume now that
the iterations of the \inneralg\ algorithm at Step~2.(a) of the outer
iteration $j$ are numbered $(j,0), \,(j,1), \ldots, (j, e_j)$ and note that
the mechanism of the \outeralg\ algorithm ensures that iteration $(j,e_j)$ is
successful for all $j$.  Now define, 
for $k \in \calK_+$, 
\beqn{Ikdef}
\calI_k
\eqdef \{ (k,0), \ldots, (k,e_k), \ldots, (j,0), \ldots, (j,e_j), \ldots, (k+\ell(k)-1,0), \ldots (k+\ell(k)-1,e_{k+\ell(k)-1}) \}
\eeqn
the index set of all inner iterations necessary to deduce $x_{k+\ell(k)}$ from
$x_k$.  Observe that, by the definitions  \req{K-pos} and \req{Ikdef}, the
index set of all inner iterations before termination is given by $\cup_{k \in
  \calK_+} \calI_k$, and therefore that the number of evaluations of problem's
functions required to terminate in Phase~2 is bounded above by  
\beqn{rough-bound}
|\bigcup_{k \in \calK_+} \calI_k| +1
\leq \Big([f_{\rm up}-f_{\rm low}+1]\ep^{-\frac{p+1}{p}} \times \max_{k \in \calK_+} | \calI_k| \Big)+1,
\eeqn
where we added 1 to take the final evaluation into account and where
we used Lemma~\ref{number-outer} to deduce the inequality.
We now invoke the complexity properties of the \inneralg\ algorithm applied on
problem \req{phase2-problem} to obtain an upper bound on the cardinality of
each $\calI_k$.   

\llem{number-inner}{
Suppose that AS.3--AS.8  hold. Then, for each $k \in \calK_+$ before
termination,  
\[
|\calI_k| \leq \kap{CC}^\mu\,\ep\ed^{-\frac{p+1}{p}}, 
\]
where $\kap{CC}^\mu$ is independent of $\ep$ and $\ed$ and captures the
problem-dependent constants associated with problem \req{phase2-problem} for
all values of $t_k$ satisfying  \req{tk-bounded}. 
}

\proof{
Observe first that \req{nreseps} implies that, for each $k$, $\calL(x_k)
\subseteq \calC_\beta$. Hence, because of Lemma~\ref{uniform-Lsig}, we may
apply the \inneralg\ algorithm for the minimization of $\mu(x,t_j)$ for each
$j$ such that $k \leq j < k+\ell(k)$,  Observe that \req{r-swap} guarantees
the decreasing nature of the sequence $\{\|r(x_j,t_j)\|\}_{j =
  k}^{k+\ell(k)-1}$ and hence of the sequence $\{\|r(x_{j,s},t_j)\|\}_{(j,s)
  \in \calI_k}$. For each $k \in \calK_+$, this minimization starts from the
initial value $\half \|r(x_{k,0},t_k)\|^2= \half\ep^2$ and is carried out for
all iterations with index in $\calI_k$ at worst down to  the value
$\half(\ep-\ep^{(p+1)/p})^2$ (see the first part of \req{ph2term}). We may
then invoke Lemmas~\ref{uniform-Lsig} and \ref{decrsucc} to deduce that, if
$(j,s)\in \calI_k$ is the index of a successful inner iteration and as long as
the third part of \req{ph2term} does not hold, 
\[
\begin{array}{ll}
\half \|r(x_{j,s},t_j)\|^2 - \half \|r(x_{j,s+1},t_j)\|^2 \geq \kap{CC}^{\mu,s} (\ep\ed)^{\frac{p+1}{p}},
&
\tim{ for } 0 \leq s<e_j, \tim{and}\\*[2ex]
\half \|r(x_{j,e_j},t_j)\|^2 - \half \|r(x_{j+1,0},t_{j+1})\|^2 \geq \kap{CC}^{\mu,s} (\ep\ed)^{\frac{p+1}{p}}& \\
\end{array}
\]
for some constant $\kap{CC}^{\mu,s} > 0$ independent of $\ep$, $\ed$, $s$ and $j$.
As a  consequence, the number of successful iterations  of the
\inneralg\ algorithm needed to compute $x_{k+\ell(k)}$ from $x_k$ cannot
exceed 
\[ 
\kap{CC}^{\mu,s}
     \left[\frac{\ep^2 - (\ep-\ep^{\frac{p+1}{p}})^2}{2(\ep\ed)^{\frac{p+1}{p}}}\right]
= \kap{CC}^{\mu,s} 
    \left[\frac{2\ep\ep^{\frac{p+1}{p}}-\ep^{2\frac{p+1}{p}}}{2\ep^{\frac{p+1}{p}}\ed^{\frac{p+1}{p}}}\right]
< \kap{CC}^{\mu,s} \ep\ed^{-\frac{p+1}{p}}.
\]
We now use Lemma~\ref{uniform-Lsig} again and invoke Lemma~\ref{unsuccbounbd}
to account for possible unsuccessful inner iterations, yielding that the total
number of successful and unsuccessful iterations of the \inneralg\ algorithm
necessary to deduce $x_{k+\ell(k)}$ from $x_k$ is bounded above by 
\[
\kappa_u \,\kap{CC}^{\mu,s}\, \ep\ed^{-\frac{p+1}{p}}
\eqdef \kap{CC}^\mu\, \ep\ed^{-\frac{p+1}{p}}.
\]
} 

\noindent
We finally combine our results in a final theorem stating our evaluation
complexity bound for the \outeralg\ algorithm applied on the general smooth
nonlinear optimization problem. 

\lthm{general-th}{
Suppose that AS.3--AS.8 hold. Then, for some constants $\kap{CC}^c$ and
$\kap{CC}^\mu$ independent of $\ep$ and $\ed$, the \outeralg\ algorithm
applied on problem \req{genprob} needs at most 
\beqn{finbound}
\left\lfloor
\Big[\kap{CC}^c \|c(x_0)\|^2+\kap{CC}^\mu (f_{\rm up}-f_{\rm low}+1)\Big] \ep^{-\frac{1}{p}}\ed^{-\frac{p+1}{p}}
\right\rfloor
\eeqn
evaluations of $f$, $c$ and their derivatives up to order $p$ to compute a
point $x_\epsilon$ such that either 
\beqn{xinfeas}
\|c(x_\epsilon)\| >\ep -\ep^{\frac{p+1}{p}}
\tim{ and } 
\chi_{\half\|c\|^2}(x_\epsilon) \leq \ed \|c(x_\epsilon)\|
\eeqn
or
\beqn{stat}
\|c(x_\epsilon)\| \leq \ep
\tim{ and }
\chi_\Lambda(x_\epsilon,y_\epsilon)  \leq \delta \ed \| (y_\epsilon, 1)\|,
\eeqn
where $\Lambda(x,y) \eqdef f(x) + y^Tc(x)$ is the Lagrangian with respect
to the equality constraints and $y_\epsilon$ is a vector of Lagrange
multipliers associated with the equality constraints. 
}

\proof{
If the \outeralg\ algorithm terminates in Phase~1, we immediately obtain that
\req{xinfeas} holds, and Lemma~\ref{ph1-compl} then ensures that the number of
evaluations of $c$ and its derivatives cannot exceed 
\beqn{bph1}
\kap{CC}^c \|c(x_0)\| \, \ep^{-\frac{1}{p}} \, \ed^{-\frac{p+1}{p}}.
\eeqn
The conclusions of the theorem therefore hold in this case.

Let us now assume that termination does not occur in Phase~1.  Then the
\outeralg\ algorithm must terminate after a number of evaluations of $f$ and
$c$ and their derivatives which is bounded above by the upper bound on the
number of evaluations in Phase~1 given by \req{bph1} plus the bound on the
number of evaluations of $\mu$ given by \req{rough-bound} and
Lemma~\ref{number-inner}. This yields the combined upper bound
\[
\kap{CC}^c \|c(x_0)\|^2\,\ep^{-\frac{1}{p}} \, \ed^{-\frac{p+1}{p}}
+ \left[\kap{CC}^\mu \ep \,\ed^{-\frac{p+1}{p}}\right]\,
  \left[ (f_{\rm up}-f_{\rm low}+1) \,\ep^{-\frac{p+1}{p}}\right],
\]
and  \req{finbound} follows. 
Remember now that \req{conds-at-termination} holds at termination of Phase~2,
and therefore that 
\beqn{rbiggish}
\ep \geq \|r(x_\epsilon,t_\epsilon)\| 
\geq  \ep-\ep^{\frac{p+1}{p}}
= \ep \left(1 - \epsilon^{\frac{1}{p}}\right)
\geq \delta^{-1} \ep.
\eeqn
Moreover,  we also obtain from \req{conds-at-termination} that
\beqn{chi-small}
\chi_\mu(x_\epsilon,t_\epsilon) 
\leq \ep\ed
\leq \delta \ed \|r(x_\epsilon,t_\epsilon)\|.
\eeqn
Assume first that $f(x_\epsilon) = t_\epsilon$.  Then, using \req{chidef} and
the definition of $r(x,t)$,, we deduce that 
\[
\chi_{\half\|c\|^2}(x_\epsilon) 
= \chi_\mu(x_\epsilon,t_k) 
\leq\delta \ed \|c(x_\epsilon)\|
\]
and \req{xinfeas} is again satisfied.  Assume now that $f(x_\epsilon) >
t_\epsilon$ (the case where $f(x_\epsilon) < t_\epsilon$ is excluded by
\req{conds-at-termination}).  Defining now
\[
y_\epsilon \eqdef \frac{c(x_\epsilon)}{f(x_\epsilon)-t_\epsilon}.
\]
and successively using the definition of $\Lambda(x,y)$, the inequality
$f(x_\epsilon) \geq t_\epsilon$, the linearity of $\chi_\Lambda$ for positive
multiples of $\nabla_x \Lambda(x,y)$ (see \req{chidef}), \req{chi-small} and
the definition of $r(x,t)$, we deduce that
\[
\chi_\Lambda(x_\epsilon,y_\epsilon) 
=  \bigfrac{\chi_\mu(x_\epsilon,t_\epsilon)}{f(x_\epsilon)-t_\epsilon} 
\leq  \delta \ed \bigfrac{\|r(x_\epsilon,t_\epsilon)\|}{f(x_\epsilon)-t_\epsilon}
 =    \delta \ed \| (y_\epsilon, 1)\|
\]
This finally implies  \req{stat} since
$\|c(x_\epsilon)\| \leq \|r(x_\epsilon,t_\epsilon)\| \leq \ep$.
} 

Note that the bound \req{finbound} is $O(\epsilon^{-\frac{p+2}{p}})$ whenever
$\ep = \ed = \epsilon$. Note also that, because of \req{rbiggish}, the first
inequality in \req{xinfeas} ensures that $\|c(x_\epsilon)\| \geq \delta^{-1}
\ep$.

\numsection{Discussion}\label{concl-s}

We have first shown in Section~\ref{convex-s} that, if derivatives of the
objective function up to order $p$ can be evaluated and if the $p$-th one is
Lipschitz continuous, then the \inneralg\ algorithm applied of the convexly
constrained problem \req{prob} needs at most $O(\epsilon^{\frac{p+1}{p}})$
evaluations of $f$ and its derivatives to compute an $\epsilon$-approximate
first-order critical point. This worst-case bound corresponds to that obtained
in \cite{CartGoulToin12b} when $p = 2$, but with significantly weaker
assumptions.  Indeed, the present proposal no longer needs any assumption on
the number of descent steps in the subproblem solution, the iterates are no
longer assumed to remain in a bounded set and the Lipschitz continuity of the
gradient is no longer necessary.  That these stronger results are obtained as
the result of a considerably simpler analysis is an added bonus. While we have
not developed here the case (covered for $p=2$ in \cite{CartGoulToin12b})
where the $p$-th derivative is only known approximately (in the sense that
$\nabla_x^p f(x_k)$ is replaced in the model's expression by some tensor $B_k$
such that the norm of $(\nabla_x^p f(x_k) - B_k)$ applied $p-1$ times to $s_k$
must be $O(\|s_k\|^p)$), the generalization of the present proposal to cover
this situation is easy.

The proposed worst-case evaluation bound also generalizes that of
\cite{BirgGardMartSantToin17} for unconstrained optimization to the
case of set-constrained problems, under very weak assumptions on the
feasible set. As was already the case for $p\leq 2$, it is remarkable
that the complexity bound for the considered class of problems (which
includes the standard bound constrained case) is, for all $p \geq
1$, identical in order to that of unconstrained problems.

The present framework for handling convex constraints is however not free of
limitations, resulting from the choice to transfer difficulties associated
with the original problem to the subproblem solution, thereby sparing precious
evaluations of $f$ and its derivatives.  The first is that we need to compute
values of $\chi_f$ and $\chi_{m_k}$.  While this is straightforward for simple
convex sets such boxes, the process might be more intensive for the general
case, although the $\|\cdot\|_\chi$ norm may be chosen to simplify this
computation.  The second limitation is that the approximate solution of the
subproblem may also be very expensive in terms of internal calculations (we do
not consider here suitable algorithms for this purpose).  Observe nevertheless
that, crucially, neither the computation of the criticality measures nor the
subproblem solution involve evaluating the objective function or its
derivatives: despite their potential computational drawbacks, they have
therefore no impact on the evaluation complexity of the original
problem. Moreover, as the cost of evaluating any constraint
function/derivative possibly necessary for computing $\chi_f$ and $\chi_{m_k}$
is neglected by the present approach, it must therefore be seen as a suitable
framework to handle "cheap inequality constraints" such as simple bounds.

We have also shown in Section~\ref{general-s} that the evaluation complexity
of finding an approximate first-order scaled critical point for the general
smooth nonlinear optimization problem involving both equality and inequality
constraints is at most $O(\ep^{-1/p}\ed^{-(p+1)/p})$ evaluations of the
objective function, constraints and their derivatives up to order $p$. We
refer here to an "approximate scaled critical point" in that such a point is
required to satisfy \req{xinfeas} or \req{stat}, where the accuracy is scaled
by the size of the constraint violation or that of the Lagrange multipliers.
Because this bound now only depends on the assumptions necessary to prove the
evaluation complexity bound for the \inneralg\ algorithm in
Section~\ref{convex-s}, it therefore strengthens and generalizes that of
\cite{CartGoulToin13a} since the latter directly hinges on
\cite{CartGoulToin12b}.

Interestingly, an $O(\ep \ed^{-(p+1)/p} \min[\ed,\ep]^{-(p+1)/p})$ evaluation
complexity bound was also proved by Birgin, Gardenghi, Mart\'{i}nez, Santos
and Toint in \cite{BirgGardMartSantToin16} for \emph{unscaled}, standard KKT
conditions and in the least expensive of three cases depending on the degree
of degeneracy identifiable by the algorithm. Even if the bounds for the scaled
and unscaled cases coincide in order when $\ep \leq \ed$, comparing the two
results is however not straightforward.  On one hand the scaled conditions
take into account the possibly different scaling of the objective function and
constraints.  On the other hand the same scaled conditions may result in
earlier termination with \req{stat} if the Lagrange multipliers are very
large, as \req{stat} is then consistent with the weaker requirement of finding
a John's point.  But the framework discussed in the present paper also differs
from that of \cite{BirgGardMartSantToin16} in additional significant ways.
The first is that the present one provides a potentially stronger version of
the termination of the algorithm at infeasible points (in Phase~1): indeed the
second part of \req{xinfeas} can be interpreted as requiring that the size of
the feasible linear decrease of $\|c(x)\|$ is below $\ed$, while
\cite{BirgGardMartSantToin16} considers the gradient of $\|c(x)\|^2$
instead. The second is that, if termination occurs in Phase~2 for an
$x_\epsilon$ such that $\chi_{\|\cdot\|^2}(x_\epsilon) =
\|J(x_\epsilon)^Tc(x_\epsilon)\|$ is itself of order $\ep\ed$ (thereby
covering the case where $f(x_\epsilon)=t_k$ discussed in
Theorem~\ref{general-th}) , then Birgin \textit{et al.} show that the
{\L}ojaciewicz inequality \cite{Loja65} must fail for $c$ in the limit for
$\ep$ and $\ed$ tending to zero (see \cite{BirgGardMartSantToin16} for
details).  This observation is interesting because smooth functions satisfy
the {\L}ojaciewicz inequality under relatively weak conditions, implying that
termination in these circumstances is unlikely. The same information is also
obtained in \cite{BirgGardMartSantToin16}, albeit at the price of worsening
the evaluation complexity bound mentioned above by an order of magnitude in
$\ed$.  We also note that the approach of \cite{BirgGardMartSantToin16}
requires the minimization, at each iteration, of a residual whose second
derivatives are discontinuous, while all functions used in the present paper
are $p$ times continuously differentiable.  A final difference between the two
approaches is obviously our introduction of $\chi_\Lambda$ and $\chi_{\half
  \|c\|^2}$ in the expression of the criticality condition in
Theorem~\ref{general-th} for taking the inequality constraints into account.

We conclude by recalling that parallel results are obtained in
\cite{CartGoulToin15a} using 
\[
\pi_f(x) = \|P_{\calF}[x-\nabla_x^1f(x)] - x \|
\]
as an alternative criticality measure replacing $\chi_f(x)$.

{\footnotesize

}

\bibliographystyle{plain}
\bibliography{/home/pht/bibs/refs}
\end{document}